\documentclass{article}%
\usepackage{makeidx}
\usepackage{amsfonts}
\usepackage{amsmath}%
\usepackage{amssymb}%
\usepackage{graphicx}
\newtheorem{theorem}{Theorem}

\newenvironment{proof}[1][Proof]{\noindent\textbf{#1.} }{\ \rule{0.5em}{0.5em}}
\begin{document}

\title{Chromogeometry}
\author{N J Wildberger\\School of Mathematics and Statistics\\UNSW Sydney 2052 Australia}
\maketitle
\date{}

\begin{abstract}
Chromogeometry brings together Euclidean geometry (called \textit{blue}) and
two relativistic geometries (called \textit{red }and \textit{green}), in a
surprising three-fold symmetry. We show how the red and green `Euler lines'
and `nine-point circles' of a triangle interact with the usual blue ones, and
how the three orthocenters form an associated triangle with interesting
collinearities. This is developed in the framework of rational trigonometry
using \textit{quadrance} and \textit{spread} instead of \textit{distance} and
\textit{angle}. The former are more suitable for relativistic geometries.

\end{abstract}

\section*{Introduction}

Three-fold symmetry is at the heart of a lot of interesting mathematics and
physics. This paper shows that it also plays an unexpected role in planar
geometry, in that the familiar Euclidean geometry is only one of a trio of
interlocking metrical geometries. We refer to Euclidean geometry here as
\textit{blue }geometry; the other two geometries, called \textit{red }and
\textit{green,} are relativistic in nature and are associated with the names
of Lorentz, Einstein and Minkowski.

The three geometries support each other and interact in a rich way. This
transcends Klein's Erlangen program, since there are now \textit{three groups
}acting on a space. Remarkable algebraic identities lie at the heart of the explanations.

The results described here are just the tip of an iceberg, leading to many
rich generalizations of results of Euclidean geometry, with much waiting to be
discovered and explored, see for example \cite{Wildberger 4} for applications
to conics and \cite{Wildberger3} for connections with one dimensional metrical geometry.

The basic structure of all three geometries are the same---they are ruled by
the laws of \textit{rational trigonometry} as developed recently in
\cite{Wildberger}, which hold over a general field not of characteristic two.
Although over the rational numbers (or the `real numbers') there are
significant differences between the Euclidean (blue) version and the other two
(red and green), it is the \textit{interaction }of all three which yields the
biggest surprises.

To start the ball rolling, this paper first introduces the phenomenon in the
context of the classical Euler line and nine-point circle of a triangle. Then
we recall the main laws of rational trigonometry, introduce the basic facts
about the three geometries and state some explicit formulas, and then show how
chromogeometry allows us to enlarge our understanding of the geometry of a
triangle. In particular we associate to each triangle $\overline{A_{1}%
A_{2}A_{3}}$ in the Cartesian plane a second interesting triangle which we
call the $\Omega$\textit{-triangle }of $\overline{A_{1}A_{2}A_{3}}$.

The results are verified by routine but sometimes lengthy computation, and
they inevitably reduce to algebraic identities, some of which are quite
lovely. The development takes place in the framework of \textit{universal
geometry}, so that we are interested primarily in what happens over
\textit{arbitrary fields}. The paper (\cite{Wildberger2}) shows that universal
geometry also extends to arbitrary quadratic forms, and embraces both
spherical and hyperbolic geometries in a projective version.

\section*{Euler lines and nine-point circles in relativistic settings}

Recall that for a triangle $\overline{A_{1}A_{2}A_{3}}$ the intersection of
the medians is the \textbf{centroid }$G,$ the intersection of the altitudes is
the \textbf{orthocenter }$O$ and the intersection of the perpendicular
bisectors of the sides is the \textbf{circumcenter} $C,$ which is the center
of the circumcircle of the triangle. Remarkably, it was left to Euler to
discover that these three points are collinear, and that $G$ divides
$\overline{OC}$ in the (affine) proportion $2:1.$ Furthermore the center $N$
of the circumcircle of the triangle $\overline{M_{1}M_{2}M_{3}}$ of midpoints
of the sides of $\overline{A_{1}A_{2}A_{3}}$ (called the \textbf{nine-point
circle} of $\overline{A_{1}A_{2}A_{3}}$) also lies on the Euler line, and is
the midpoint of $\overline{OC}.$%

%TCIMACRO{\FRAME{dhF}{3.0564in}{3.3698in}{0pt}{}{}{chromoex2blue.eps}%
%width 3.0564in;  height 3.3698in;  depth 0pt;  original-width 4.0365in;
%original-height 4.4524in;  cropleft "0";  croptop "1";  cropright "1";
%cropbottom "0";  filename 'ChromoEx2Blue.eps';file-properties "XNPEU";}}}%
%BeginExpansion
\begin{center}
\includegraphics[
height=3.3698in,
width=3.0564in
]%
{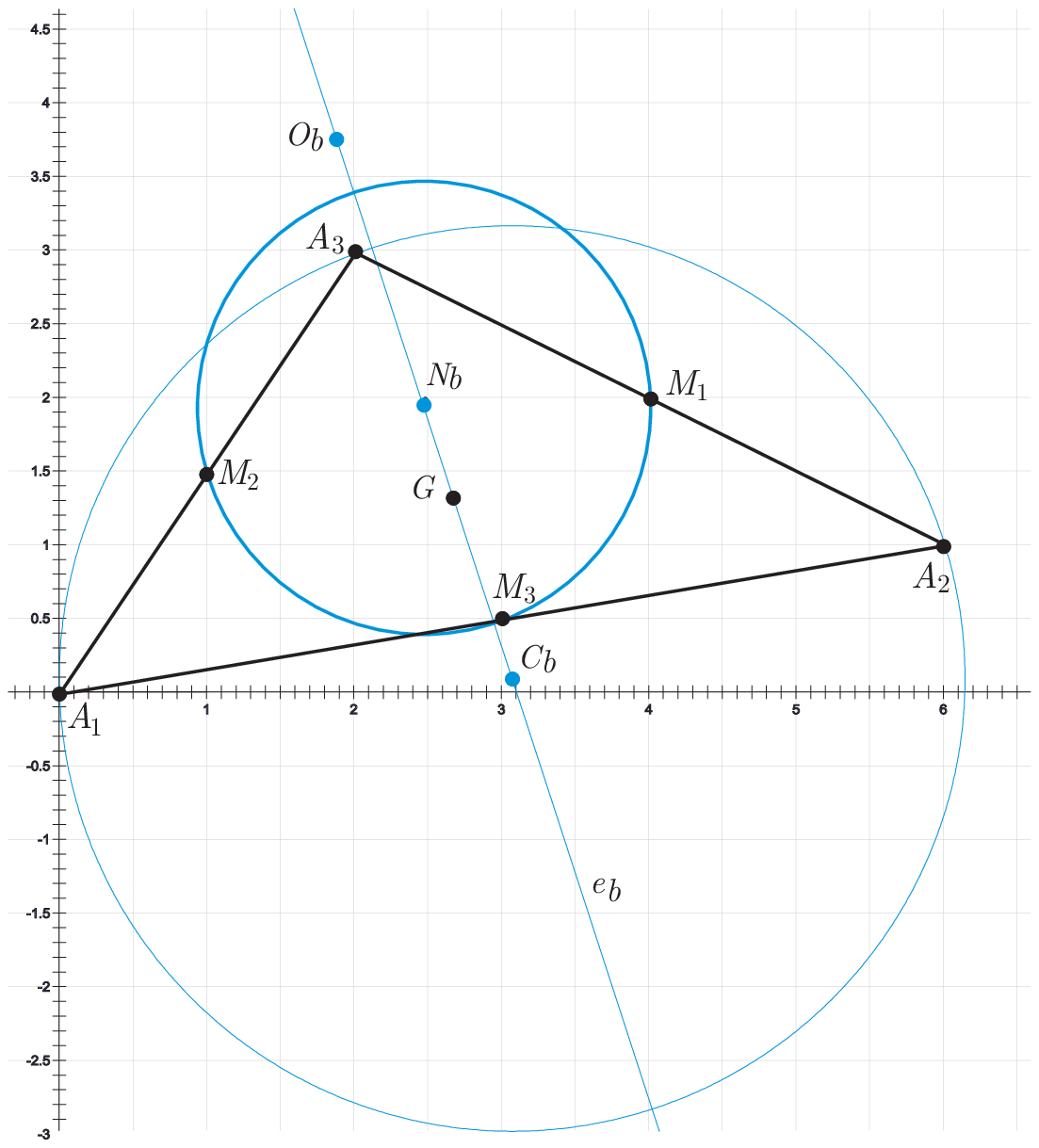}%
\end{center}
%EndExpansion
This is above shown for the triangle $\overline{A_{1}A_{2}A_{3}}$ with points%
\[%
\begin{tabular}
[c]{lllll}%
$A_{1}\equiv\left[  0,0\right]  $ &  & $A_{2}\equiv\left[  6,1\right]  $ &  &
$A_{3}\equiv\left[  2,3\right]  .$%
\end{tabular}
\]
The triangle $\overline{A_{1}A_{2}A_{3}}$ is in black, while the circumcircle
and nine-point circle are in blue (the latter more boldly), as are the Euler
line and the points $O,C$ and $N,$ which are given the subscript $b$ for blue,
and henceforth referred to as the \textbf{blue Euler line}, the \textbf{blue
orthocenter }etc.

Planar Euclidean geometry rests on the \textbf{blue quadratic form}
$x^{2}+y^{2}$ (or if you prefer the corresponding symmetric bilinear form, or
dot product). It is also interesting to consider the \textbf{red quadratic
form} $x^{2}-y^{2}$ which figures prominently in two dimensional special
relativity. In this case, two lines are red perpendicular precisely when one
can be obtained from the other by ordinary Euclidean reflection in a
\textbf{red null line, }which is red perpendicular to itself, and has usual
slope $\pm1.$

It turns out that for any triangle $\overline{A_{1}A_{2}A_{3}}$ the three red
altitudes also intersect, now in a point called the \textbf{red orthocenter}
and denoted $O_{r},$ and the three perpendicular bisectors also intersect in a
point called the \textbf{red circumcenter} and denoted $C_{r}$. This latter
point is the center of the unique red circle through the three points of the
triangle, where a \textbf{red circle} is given by an equation of the form
\[
\left(  x-x_{0}\right)  ^{2}-\left(  y-y_{0}\right)  ^{2}=K.
\]
This is what we would usually call a rectangular hyperbola, with axes in the
red null directions.%

%TCIMACRO{\FRAME{dhF}{2.9795in}{3.2352in}{0pt}{}{}{chromoex2red.eps}%
%width 2.9795in;  height 3.2352in;  depth 0pt;  original-width 4.18in;
%original-height 4.5422in;  cropleft "0";  croptop "1";  cropright "1";
%cropbottom "0";  filename 'ChromoEx2Red.eps';file-properties "XNPEU";}}}%
%BeginExpansion
\begin{center}
\includegraphics[
height=3.2352in,
width=2.9795in
]%
{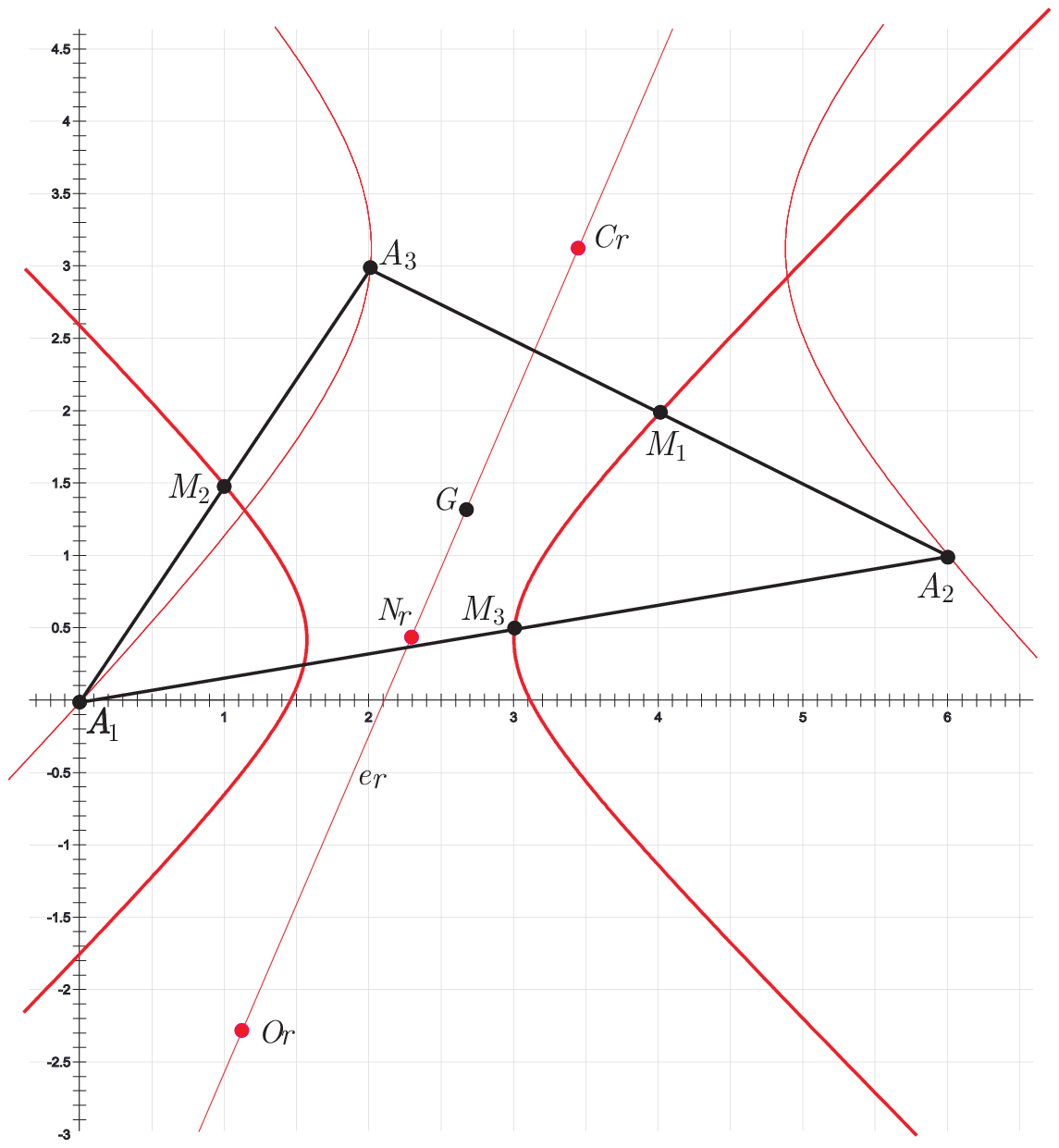}%
\end{center}
%EndExpansion
This diagram shows the same triangle $\overline{A_{1}A_{2}A_{3}}$, as well as
the red circumcircle, the red nine-point circle and the red orthocenter,
circumcenter, nine-point center and centroid $G$, the latter being independent
of colour.

Note that these points all lie on a line---the \textbf{red Euler line}, and
the affine relationships between these points is exactly the same as for the
blue Euler line, so that for example $N_{r}$ is the midpoint of $\overline
{O_{r}C_{r}}$.

In the classical framework, there are some difficulties in setting up this
relativistic geometry, as `distance' and `angle' are problematic. In universal
geometry one regards the \textit{quadratic form }as primary, not its
\textit{square root}, and by expressing everything in terms of the algebraic
concepts of \textit{quadrance }and \textit{spread}, Euclidean geometry can be
built up so as to allow generalization to the relativistic framework, and
indeed to geometries built from other quadratic forms.

This approach was introduced recently in \cite{Wildberger}, see also
\cite{Wildberger5}, and works over a general field with characteristic two
excluded for technical reasons, as shown in \cite{Wildberger2}. The
possibility of relativistic geometries over other fields seems particularly attractive.

There is a third geometry, that associated to the \textbf{green quadratic
form} $2xy.$ Two lines are green perpendicular when one is the ordinary
Euclidean reflection of the other in a line parallel to the axes, the latter
being a \textbf{green null line}.

Since $x^{2}-y^{2}$ and $2xy$ are conjugate by a simple change of variable, it
should be no surprise that the corresponding relations between the green
orthocenter $O_{g},$ green circumcenter $C_{g}$, green nine-point center
$N_{g}$ and the centroid $G$ hold as well. Here is the relevant diagram for
our triangle $\overline{A_{1}A_{2}A_{3}}.$%

%TCIMACRO{\FRAME{dhF}{3.1133in}{3.2175in}{0pt}{}{}{chromoex2green.eps}%
%width 3.1133in;  height 3.2175in;  depth 0pt;  original-width 4.3194in;
%original-height 4.466in;  cropleft "0";  croptop "1";  cropright "1";
%cropbottom "0";  filename 'ChromoEx2Green.eps';file-properties "XNPEU";}}}%
%BeginExpansion
\begin{center}
\includegraphics[
height=3.2175in,
width=3.1133in
]%
{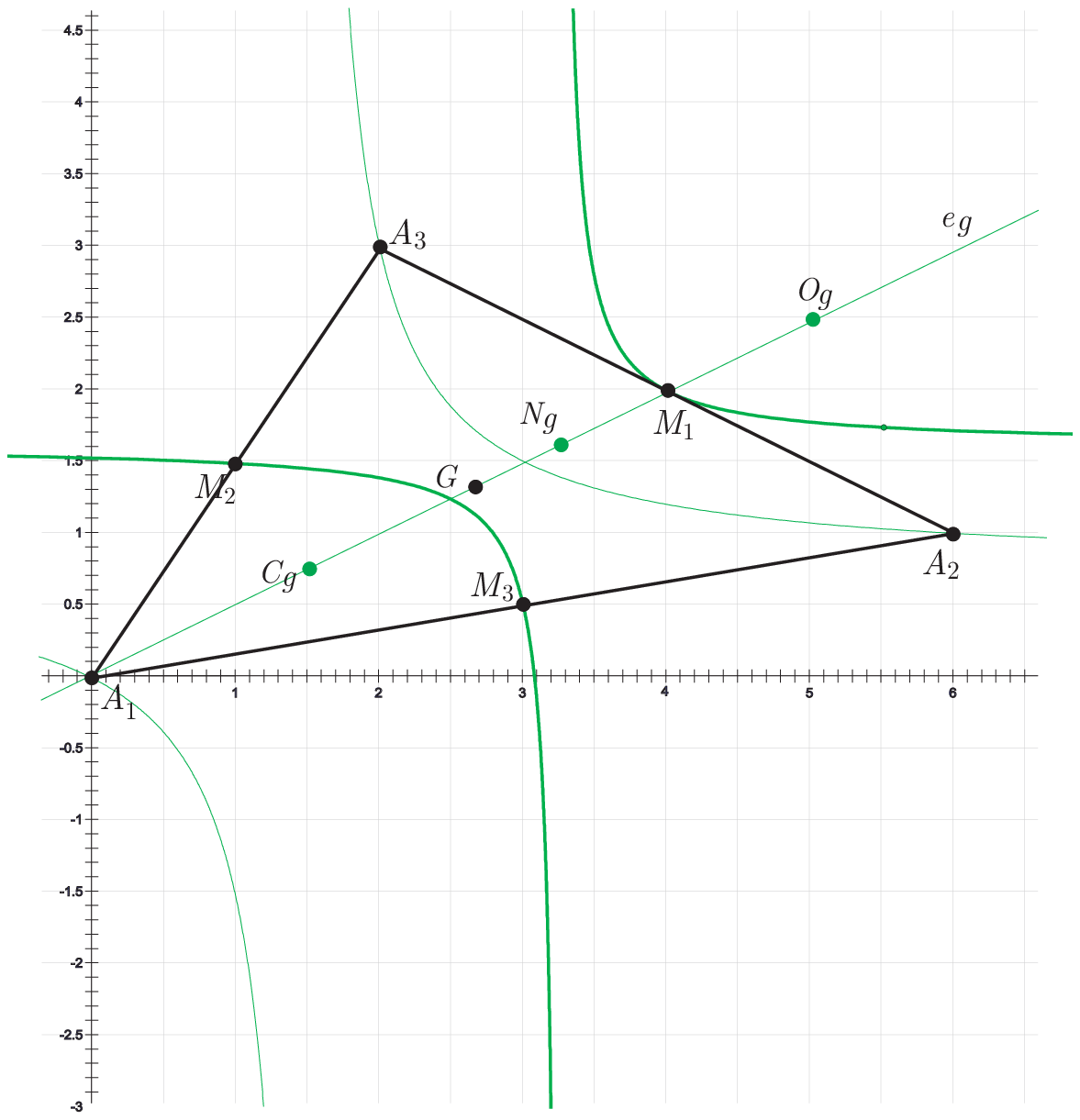}%
\end{center}
%EndExpansion

I conjecture that \textit{most theorems of planar Euclidean geometry, when
formulated algebraically in the context of universal geometry, extend to the
red and green situations}. However there are exceptions. For example, over the
`real numbers' there are no equilateral triangles in the red or green
geometries, so Napoleon's theorem and Morley's theorem will not have direct analogs.

Much could be said further to support this conjecture, but this is not what I
wish to pursue here. Instead, let's consider a completely new phenomenon.
Observe what happens when the three diagrams are put together!%

%TCIMACRO{\FRAME{dhF}{3.5373in}{3.7176in}{0pt}{}{}{chromoex2.eps}%
%width 3.5373in;  height 3.7176in;  depth 0pt;  original-width 4.3194in;
%original-height 4.5422in;  cropleft "0";  croptop "1";  cropright "1";
%cropbottom "0";  filename 'ChromoEx2.eps';file-properties "XNPEU";}}}%
%BeginExpansion
\begin{center}
\includegraphics[
height=3.7176in,
width=3.5373in
]%
{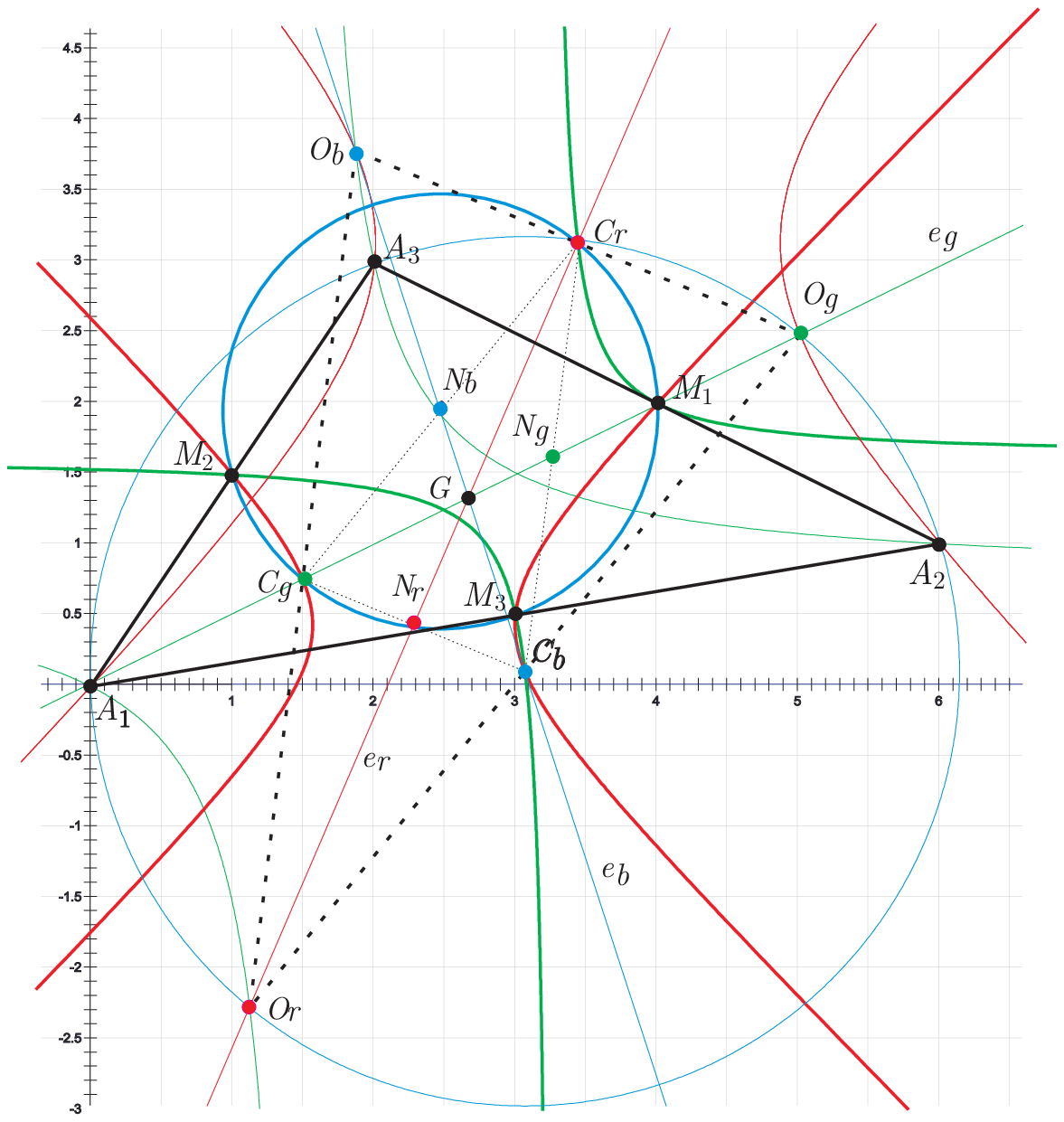}%
\end{center}
%EndExpansion

We get remarkable collinearities, for example between $O_{b},$ $C_{r}$ and
$O_{g},$ and between $C_{b},$ $N_{r}$ and $C_{g}$, with furthermore $C_{r}$
the midpoint of $\overline{O_{b}O_{g}},$ and $N_{r}$ the midpoint of
$\overline{C_{b}C_{g}}$. Also we observe that, for example, $O_{r}$ and $O_{g}
$ lie on the blue circumcircle of $\overline{A_{1}A_{2}A_{3}},$ while $C_{r}$
and $C_{g}$ lie on the blue nine-point circle of $\overline{A_{1}A_{2}A_{3}}$.
The three colours generally interact symmetrically, so the same relations hold
if we permute colours.

However there are some aspects of chromogeometry in which this symmetry is
broken. The blue geometry as we shall see behaves somewhat differently from
the red and the green in certain contexts, and when we come to explicit
formulas we will see that the green geometry is often simpler. The red
geometry seems less inclined to distinguish itself.

\section*{Rational trigonometry}

Let's now proceed more formally, beginning with the main definitions and laws
of rational trigonometry. We work over a fixed field, not of characteristic
two, whose elements will be called \textbf{numbers}. The \textbf{plane }will
consist of the standard vector space of dimension two over this field. A
\textbf{point}, or \textbf{vector}, is an ordered pair $A\equiv\left[
x,y\right]  $ of numbers. The origin is denoted $O\equiv\left[  0,0\right]  $.

A \textbf{line} is a proportion $l\equiv\left\langle a:b:c\right\rangle $
where $a$ and $b$ are not both zero. The point $A\equiv\left[  x,y\right]  $
\textbf{lies on} the line $l\equiv\left\langle a:b:c\right\rangle $, or
equivalently the line $l$ \textbf{passes through} the point $A$, precisely
when%
\[
ax+by+c=0.
\]

This is not the only possible convention, and the reader should be aware that
it is prejudiced towards the usual Euclidean (blue) geometry. For any two
points $A_{1}\equiv\left[  x_{1},y_{1}\right]  $ and $A_{2}\equiv\left[
x_{2},y_{2}\right]  $ there is a unique line $l\equiv A_{1}A_{2}$ which passes
through them both. Specifically we have
\[
A_{1}A_{2}=\left\langle y_{1}-y_{2}:x_{2}-x_{1}:x_{1}y_{2}-x_{2}%
y_{1}\right\rangle .
\]
Three points $\left[  x_{1},y_{1}\right]  $, $\left[  x_{2},y_{2}\right]  $
and $\left[  x_{3},y_{3}\right]  $ are \textbf{collinear} precisely when they
lie on the same line, which amounts to the condition
\begin{equation}
x_{1}y_{2}-x_{1}y_{3}+x_{2}y_{3}-x_{3}y_{2}+x_{3}y_{1}-x_{2}y_{1}%
=0.\label{AreaFormula}%
\end{equation}
Three lines $\left\langle a_{1}:b_{1}:c_{1}\right\rangle $, $\left\langle
a_{2}:b_{2}:c_{2}\right\rangle $ and $\left\langle a_{3}:b_{3}:c_{3}%
\right\rangle $ are \textbf{concurrent }precisely when they pass through the
same point, which amounts to the condition
\[
a_{1}b_{2}c_{3}-a_{1}b_{3}c_{2}+a_{2}b_{3}c_{1}-a_{3}b_{2}c_{1}+a_{3}%
b_{1}c_{2}-a_{2}b_{1}c_{3}=0.
\]

Fix a symmetric bilinear form, denoted by the dot product $A_{1}\cdot A_{2}$.
In practise we will take this bilinear form to be non-degenerate. The line
$A_{1}A_{2}$ is \textbf{perpendicular} to the line $B_{1}B_{2}$ precisely
when
\[
\left(  A_{2}-A_{1}\right)  \cdot\left(  B_{2}-B_{1}\right)  =0.
\]
A point $A$ is a\textbf{\ null point }or\textbf{\ null vector} precisely when
$A\cdot A=0.$ The \textbf{origin} $O$ is always a null point, but there may be
others. A line $A_{1}A_{2}$ is a \textbf{null line} precisely when the vector
$A_{2}-A_{1}$ is a null vector.

A set $\left\{  A_{1},A_{2},A_{3}\right\}  $ of three distinct non-collinear
points is a \textbf{triangle} and is denoted $\overline{A_{1}A_{2}A_{3}}$. The
\textbf{lines} of the triangle are $l_{3}\equiv A_{1}A_{2},$ $l_{2}\equiv
A_{1}A_{3}$ and $l_{1}\equiv A_{2}A_{3}.$ A triangle is \textbf{non-null}
precisely when each of its lines is non-null. A \textbf{side} of the triangle
is a subset of $\left\{  A_{1},A_{2},A_{3}\right\}  $ with two elements, and
is denoted $\overline{A_{1}A_{2}}$ etc. A \textbf{vertex} of the triangle is a
subset of $\left\{  l_{1},l_{2},l_{3}\right\}  $ with two elements, and is
denoted $\overline{l_{1}l_{2}}$ etc.

The \textbf{quadrance} between the points $A_{1}$ and $A_{2}$ is the number
\[
Q\left(  A_{1},A_{2}\right)  \equiv\left(  A_{2}-A_{1}\right)  \cdot\left(
A_{2}-A_{1}\right)  .
\]
The line $A_{1}A_{2}$ is a null line precisely when $Q\left(  A_{1}%
,A_{2}\right)  =0$.

The \textbf{spread} between the non-null lines $A_{1}A_{2}$ and $B_{1}B_{2}$
is the number%
\[
s\left(  A_{1}A_{2},B_{1}B_{2}\right)  \equiv1-\frac{\left(  \left(
A_{2}-A_{1}\right)  \cdot\left(  B_{2}-B_{1}\right)  \right)  ^{2}}{Q\left(
A_{1},A_{2}\right)  Q\left(  B_{1},B_{2}\right)  }.
\]
This is independent of the choice of points lying on the two lines. Two
non-null lines are perpendicular precisely when the spread between them is
$1.$

Here are the five main laws of planar rational trigonometry in this general
setting, replacing the usual Sine law, Cosine law etc. Proofs can be found in
\cite{Wildberger2}. Aside from giving new directions to geometry, these laws
have the potential to change the teaching of high school mathematics, because
they are simpler, and allow faster and more accurate calculations in practical
problems. But the advantage for us here is that they \textit{hold for general
quadratic forms}, and in particular for each of the blue, red and green geometries.

\begin{theorem}
[Triple quad formula]The points $A_{1},A_{2}$ and $A_{3}$ are collinear
precisely when the quadrances $Q_{1}\equiv Q\left(  A_{2},A_{3}\right)  $,
$Q_{2}\equiv Q\left(  A_{1},A_{3}\right)  $ and $Q_{3}\equiv Q\left(
A_{1},A_{2}\right)  $ satisfy%
\[
\left(  Q_{1}+Q_{2}+Q_{3}\right)  ^{2}=2\left(  Q_{1}^{2}+Q_{2}^{2}+Q_{3}%
^{2}\right)  .
\]

\end{theorem}

\begin{theorem}
[Pythagoras' theorem]For $A_{1},$ $A_{2}$ and $A_{3}$ three distinct points,
$A_{1}A_{3}$ is perpendicular to $A_{2}A_{3}$ precisely when the quadrances
$Q_{1}\equiv Q\left(  A_{2},A_{3}\right)  $, $Q_{2}\equiv Q\left(  A_{1}%
,A_{3}\right)  $ and $Q_{3}\equiv Q\left(  A_{1},A_{2}\right)  $ satisfy%
\[
Q_{1}+Q_{2}=Q_{3}.
\]

\end{theorem}

\begin{theorem}
[Spread law]Suppose the non-null triangle $\overline{A_{1}A_{2}A_{3}}$ has
quadrances $Q_{1}\equiv Q\left(  A_{2},A_{3}\right)  $, $Q_{2}\equiv Q\left(
A_{1},A_{3}\right)  $ and $Q_{3}\equiv Q\left(  A_{1},A_{2}\right)  $, and
spreads $s_{1}\equiv s\left(  A_{1}A_{2},A_{1}A_{3}\right)  $, $s_{2}\equiv
s\left(  A_{2}A_{1},A_{2}A_{3}\right)  $ and $s_{3}\equiv s\left(  A_{3}%
A_{1},A_{3}A_{2}\right)  $. Then%
\[
\frac{s_{1}}{Q_{1}}=\frac{s_{2}}{Q_{2}}=\frac{s_{3}}{Q_{3}}.
\]

\end{theorem}

\begin{theorem}
[Cross law]Suppose the non-null triangle $\overline{A_{1}A_{2}A_{3}}$ has
quadrances $Q_{1}\equiv Q\left(  A_{2},A_{3}\right)  $, $Q_{2}\equiv Q\left(
A_{1},A_{3}\right)  $ and $Q_{3}\equiv Q\left(  A_{1},A_{2}\right)  $, and
spreads $s_{1}\equiv s\left(  A_{1}A_{2},A_{1}A_{3}\right)  $, $s_{2}\equiv
s\left(  A_{2}A_{1},A_{2}A_{3}\right)  $ and $s_{3}\equiv s\left(  A_{3}%
A_{1},A_{3}A_{2}\right)  $. Then
\[
\left(  Q_{1}+Q_{2}-Q_{3}\right)  ^{2}=4Q_{1}Q_{2}\left(  1-s_{3}\right)  .
\]

\end{theorem}

Note that the Cross law includes as special cases both the Triple quad formula
and Pythagoras' theorem. The next result is the algebraic analog to the sum of
the angles in a triangle formula.

\begin{theorem}
[Triple spread formula]Suppose the non-null triangle $\overline{A_{1}%
A_{2}A_{3}}$ has spreads $s_{1}\equiv s\left(  A_{1}A_{2},A_{1}A_{3}\right)
$, $s_{2}\equiv s\left(  A_{2}A_{1},A_{2}A_{3}\right)  $ and $s_{3}\equiv
s\left(  A_{3}A_{1},A_{3}A_{2}\right)  $. Then%
\[
\left(  s_{1}+s_{2}+s_{3}\right)  ^{2}=2\left(  s_{1}^{2}+s_{2}^{2}+s_{3}%
^{2}\right)  +4s_{1}s_{2}s_{3}.
\]

\end{theorem}

A useful observation is that the Triple spread formula shows that $s_{3}=1$
implies that%
\[
s_{1}+s_{2}=1.
\]

\section*{Three fold symmetry}

The vectors $A_{1}\equiv\left[  x_{1},y_{1}\right]  $ and $A_{2}\equiv\left[
x_{2},y_{2}\right]  $ are \textbf{parallel} precisely when
\[
x_{1}y_{2}-x_{2}y_{1}=0.
\]
We will be interested in three main examples of symmetric bilinear forms.
Define the \textbf{blue} dot product
\[
\left[  x_{1},y_{1}\right]  \cdot_{b}\left[  x_{2},y_{2}\right]  \equiv
x_{1}x_{2}+y_{1}y_{2},
\]
the \textbf{red} dot product%
\[
\left[  x_{1},y_{1}\right]  \cdot_{r}\left[  x_{2},y_{2}\right]  \equiv
x_{1}x_{2}-y_{1}y_{2}%
\]
and the \textbf{green} dot product%
\[
\left[  x_{1},y_{1}\right]  \cdot_{g}\left[  x_{2},y_{2}\right]  =x_{1}%
y_{2}+x_{2}y_{1}.
\]

Note that between them these four expressions give all possible bilinear
expressions in the two vectors that involve only coefficients $\pm1,$ up to
sign. Two lines $l_{1}$ and $l_{2}$ are \textbf{blue},\textbf{\ red} and
\textbf{green perpendicular} respectively precisely when they are
perpendicular with respect to the blue, red and green forms. For lines
$l_{1}\equiv\left\langle a_{1}:b_{1}:c_{1}\right\rangle $ and $l_{2}%
\equiv\left\langle a_{2}:b_{2}:c_{2}\right\rangle $ these conditions amount to
the respective conditions
\[
a_{1}a_{2}+b_{1}b_{2}=0~~\mathrm{[blue]}%
\]
\textbf{\ }
\[
a_{1}a_{2}-b_{1}b_{2}=0~~\mathrm{[red]}%
\]
and
\[
a_{1}b_{2}+b_{1}a_{2}=0~~\mathrm{[green].}%
\]

In terms of coordinates, the formulas for the\textbf{\ blue}, \textbf{red} and
\textbf{green quadrances} between points $A_{1}\equiv\left[  x_{1}%
,y_{1}\right]  $ and $A_{2}\equiv\left[  x_{2},y_{2}\right]  $ are%
\begin{align*}
Q_{b}\left(  A_{1},A_{2}\right)   & =\left(  x_{2}-x_{1}\right)  ^{2}+\left(
y_{2}-y_{1}\right)  ^{2}\\
Q_{r}\left(  A_{1},A_{2}\right)   & =\left(  x_{2}-x_{1}\right)  ^{2}-\left(
y_{2}-y_{1}\right)  ^{2}\\
Q_{g}\left(  A_{1},A_{2}\right)   & =2\left(  x_{2}-x_{1}\right)  \left(
y_{2}-y_{1}\right)  .
\end{align*}

\begin{theorem}
[Coloured quadrances]For any points $A_{1}$ and $A_{2}$ let $Q_{b}$, $Q_{r}$
and $Q_{g}$ be the blue, red and green quadrances between $A_{1}$ and $A_{2}$
respectively. Then%
\[
Q_{b}^{2}=Q_{r}^{2}+Q_{g}^{2}.
\]

\end{theorem}

\begin{proof}
This is a consequence of the identity%
\[
\left(  r^{2}+s^{2}\right)  ^{2}=\left(  r^{2}-s^{2}\right)  ^{2}+\left(
2rs\right)  ^{2}.
\]

\end{proof}

The formulas for the \textbf{blue}, \textbf{red} and \textbf{green spreads}
between lines $l_{1}\equiv\left\langle a_{1}:b_{1}:c_{1}\right\rangle $ and
$l_{2}\equiv\left\langle a_{2}:b_{2}:c_{2}\right\rangle $ are%
\begin{align*}
s_{b}\left(  l_{1},l_{2}\right)   & =1-\frac{\left(  a_{1}a_{2}+b_{1}%
b_{2}\right)  ^{2}}{\left(  a_{1}^{2}+b_{1}^{2}\right)  \left(  a_{2}%
^{2}+b_{2}^{2}\right)  }=\allowbreak\frac{\left(  a_{1}b_{2}-a_{2}%
b_{1}\right)  ^{2}\allowbreak}{\left(  a_{2}^{2}+b_{2}^{2}\right)  \left(
a_{1}^{2}+b_{1}^{2}\right)  }\\
s_{r}\left(  l_{1},l_{2}\right)   & =1-\frac{\left(  b_{1}b_{2}-a_{1}%
a_{2}\right)  ^{2}}{\left(  b_{1}^{2}-a_{1}^{2}\right)  \left(  b_{2}%
^{2}-a_{2}^{2}\right)  }=-\frac{\left(  a_{1}b_{2}-a_{2}b_{1}\right)
^{2}\allowbreak}{\left(  a_{2}^{2}-b_{2}^{2}\right)  \left(  a_{1}^{2}%
-b_{1}^{2}\right)  }\\
s_{g}\left(  l_{1},l_{2}\right)   & =1-\frac{\left(  -a_{1}b_{2}-a_{2}%
b_{1}\right)  ^{2}}{4a_{1}a_{2}b_{1}b_{2}}=\allowbreak-\frac{\left(
a_{1}b_{2}-a_{2}b_{1}\right)  ^{2}\allowbreak}{4a_{1}a_{2}b_{1}b_{2}}.
\end{align*}
Note carefully the \textit{minus signs} that precede the final expressions in
the red and green cases.

\begin{theorem}
[Coloured spreads]For any lines $l_{1}$ and $l_{2}$ let $s_{b}$, $s_{r}$ and
$s_{g}$ be the blue, red and green spreads between $l_{1}$ and $l_{2}$
respectively. Then%
\[
\frac{1}{s_{b}}+\frac{1}{s_{r}}+\frac{1}{s_{g}}=2.
\]

\end{theorem}

\textbf{Proof.} This is a consequence of the identity%
\[
\left(  a_{1}^{2}+b_{1}^{2}\right)  \left(  a_{2}^{2}+b_{2}^{2}\right)
-\left(  a_{1}^{2}-b_{1}^{2}\right)  \left(  a_{2}^{2}-b_{2}^{2}\right)
-4a_{1}a_{2}b_{1}b_{2}=2\left(  a_{1}b_{2}-a_{2}b_{1}\right)  ^{2}.%
%TCIMACRO{\TeXButton{Proof box}{{\hspace{.1in} \rule{0.5em}{0.5em}}}}%
%BeginExpansion
{\hspace{.1in} \rule{0.5em}{0.5em}}%
%EndExpansion
\]

\section*{Quadrances}

The most important single quantity associated to a triangle $\overline
{A_{1}A_{2}A_{3}}$ with quadrances $Q_{1},Q_{2}$ and $Q_{3}$ is the
\textbf{quadrea} $\mathcal{A}$ defined by%
\[
\mathcal{A}\equiv\left(  Q_{1}+Q_{2}+Q_{3}\right)  ^{2}-2\left(  Q_{1}%
^{2}+Q_{2}^{2}+Q_{3}^{2}\right)  .
\]
By the Triple quad formula this is a measure of the non-collinearity of the
points $A_{1},A_{2}$ and $A_{3}$. We denote by $\mathcal{A}_{b},\mathcal{A}%
_{r}$ and $\mathcal{A}_{g}$ the respective \textbf{blue},\textbf{\ red
}and\textbf{\ green quadreas} of a triangle $\overline{A_{1}A_{2}A_{3}}.$

\begin{theorem}
[Quadrea ]For three points $A_{1}\equiv\left[  x_{1},y_{1}\right]  $,
$A_{2}\equiv\left[  x_{2},y_{2}\right]  $ and $A_{3}\equiv\left[  x_{3}%
,y_{3}\right]  $, the three quadreas $\mathcal{A}_{b},\mathcal{A}_{r}$ and
$\mathcal{A}_{g}$ satisfy
\[
\mathcal{A}_{b}=-\mathcal{A}_{r}=-\mathcal{A}_{g}=4\left(  x_{1}y_{2}%
-x_{1}y_{3}+x_{2}y_{3}-x_{3}y_{2}+x_{3}y_{1}-x_{2}y_{1}\right)  ^{2}.
\]

\end{theorem}

\begin{proof}
A calculation.
\end{proof}

So each quadrea of a triangle is $\pm16$ times the square of its signed area,
the latter being defined purely in an affine setting, without any need for
metrical choices.

We now adopt the convention that if no proof is given, `a calculation' is to
be assumed.

\section*{Altitudes}

\begin{theorem}
[Altitudes to a line]For any point $A$ and any line $l,$ there exist unique
lines $n_{b},n_{r}$ and $n_{g}$ through $A$ which are respectively blue, red
and green perpendicular to $l$. If $A\equiv\left[  x_{0},y_{0}\right]  $ and
$l\equiv\left\langle a:b:c\right\rangle $ then
\begin{align*}
n_{b}  & =\left\langle b:-a:-bx_{0}+ay_{0}\right\rangle \\
n_{r}  & =\left\langle b:a:-bx_{0}-ay_{0}\right\rangle \\
n_{g}  & =\left\langle a:-b:-ax_{0}+by_{0}\right\rangle .%
%TCIMACRO{\TeXButton{Proof box}{{\hspace{.1in} \rule{0.5em}{0.5em}}}}%
%BeginExpansion
{\hspace{.1in} \rule{0.5em}{0.5em}}%
%EndExpansion
\end{align*}

\end{theorem}

The lines $n_{b},n_{r},n_{g}$ are respectively the \textbf{blue},\textbf{\ red
}and\textbf{\ green altitudes} from $A$ to $l$, and they intersect $l$ at the
\textbf{feet, }provided that $l$ is non-null.

\begin{theorem}
[Perpendicularity of altitudes]For any point $A$ and any line $l,$ let
$n_{b},n_{r},n_{g}$ be the blue, red and green altitudes from $A$ to $l$
respectively. Then $n_{b}$ and $n_{r}$ are green perpendicular, $n_{r}$ and
$n_{g}$ are blue perpendicular, and $n_{g}$ and $n_{b}$ are red perpendicular.
$%
%TCIMACRO{\TeXButton{Proof box}{{\hspace{.1in} \rule{0.5em}{0.5em}}}}%
%BeginExpansion
{\hspace{.1in} \rule{0.5em}{0.5em}}%
%EndExpansion
$%
%TCIMACRO{\FRAME{dhF}{2.3248in}{1.6212in}{0pt}{}{}{altitudes.eps}%
%width 2.3248in;  height 1.6212in;  depth 0pt;  original-width 2.6261in;
%original-height 1.8223in;  cropleft "0";  croptop "1";  cropright "1";
%cropbottom "0";  filename 'Altitudes.eps';file-properties "XNPEU";}}}%
%BeginExpansion
\begin{center}
\includegraphics[
height=1.6212in,
width=2.3248in
]%
{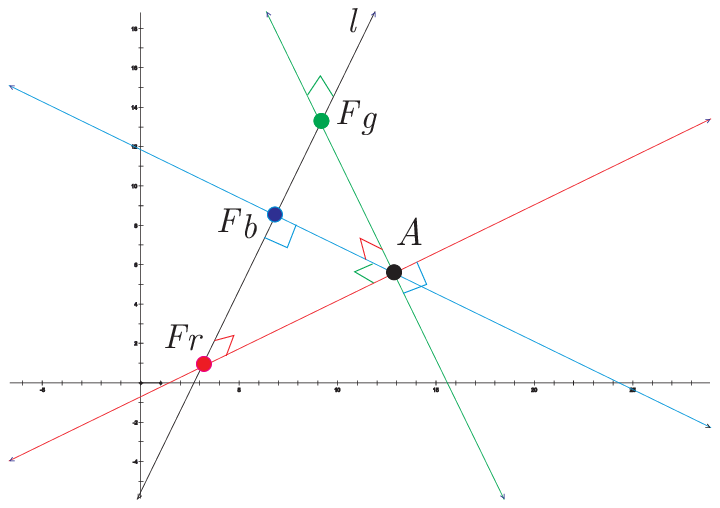}%
\end{center}
%EndExpansion

\end{theorem}

The figure shows an example of the three colour altitudes from a point $A$ to
a line $l,$ and their feet $F_{b},$ $F_{r}$ and $F_{g}.$

\begin{theorem}
[Pythagorean means]Let $l\equiv\left\langle a:b:c\right\rangle $ be a line
which is non-null in each of the three geometries. If $A$ is a point and
$F_{b},F_{r}$ and $F_{g}$ are the respective feet of the altitudes
$n_{b},n_{r}$ and $n_{g}$ from $A$ to $l$, then we have the affine relation
\[
F_{b}=\frac{\left(  a^{2}-b^{2}\right)  ^{2}}{\left(  a^{2}+b^{2}\right)
^{2}}F_{r}+\frac{4b^{2}a^{2}}{\left(  a^{2}+b^{2}\right)  ^{2}}F_{g}.
\]

\end{theorem}

\begin{proof}
Suppose that $A\equiv\left[  x_{0},y_{0}\right]  $ and $l\equiv\left\langle
a:b:c\right\rangle $. Elimination yields%
\begin{align*}
F_{b}  & =\left[  \frac{b^{2}x_{0}-aby_{0}-ac}{a^{2}+b^{2}},\frac
{-abx_{0}+a^{2}y_{0}-bc}{a^{2}+b^{2}}\right] \\
F_{r}  & =\left[  \frac{-b^{2}x_{0}-aby_{0}-ca}{a^{2}-b^{2}},\frac
{abx_{0}+a^{2}y_{0}+bc}{a^{2}-b^{2}}\right] \\
F_{g}  & =\left[  \frac{ax_{0}-by_{0}-c}{2a},\frac{-ax_{0}+by_{0}-c}%
{2b}\right]
\end{align*}
from which we deduce the result.
\end{proof}

Note again the connection with Pythagorean triples.

\section*{Anti-symmetric polynomials}

We use notation for anti-symmetric polynomials introduced in \cite[page
28]{Wildberger}. If $m$ is a monomial in the variables $x_{1},x_{2}%
,x_{3},y_{1},y_{2},y_{3},z_{1},z_{2},z_{3},\cdots$ with all indices in the
range $1,2$ and $3$, then we define $\left[  m\right]  ^{-}$ to be the
\textit{anti-symmetric }polynomial consisting of the sum of all monomials
obtained from $m$ by performing all six permutations of the indices and
multiplying each term by the sign of the corresponding permutation.

We often write such polynomials in the order described by the successive
transpositions%
\[
\left(  23\right)  ,\left(  12\right)  ,\left(  23\right)  ,\left(  12\right)
,\left(  23\right)  .
\]

For example
\begin{align*}
\left[  x_{1}y_{2}\right]  ^{-}  & \equiv x_{1}y_{2}-x_{1}y_{3}+x_{2}%
y_{3}-x_{3}y_{2}+x_{3}y_{1}-x_{2}y_{1}\\
\left[  x_{1}^{2}x_{2}y_{2}\right]  ^{-}  & \equiv x_{1}^{2}x_{2}y_{2}%
-x_{1}^{2}x_{3}y_{3}+x_{2}^{2}x_{3}y_{3}-x_{3}^{2}x_{2}y_{2}+x_{3}^{2}%
x_{1}y_{1}-x_{2}^{2}x_{1}y_{1}\\
\left[  x_{1}^{3}y_{1}\right]  ^{-}  & \equiv x_{1}^{3}y_{1}-x_{1}^{3}%
y_{1}+x_{2}^{3}y_{2}-x_{3}^{3}y_{3}+x_{3}^{3}y_{3}-x_{2}^{3}y_{2}=0.
\end{align*}
The polynomial $\left[  x_{1}y_{2}\right]  ^{-}$ is of particular importance,
since it occurs in (\ref{AreaFormula}), is twice the signed area of the
triangle $\overline{A_{1}A_{2}A_{3}}$, appears in the Quadrea theorem, and is
often a denominator in formulas in the subject.

\section*{Orthocenters}

Given a triangle $\overline{A_{1}A_{2}A_{3}},$ for each point $A_{m}$, $m=1,2
$ and $3$ we may construct the blue, red and green altitudes $a_{m}^{b}%
,a_{m}^{r}$ and $a_{m}^{g}$ respectively to the opposite side.

\begin{theorem}
[Orthocenter formulas]The three blue altitudes $a_{1}^{b},a_{2}^{b}$ and
$a_{3}^{b}$ meet in a point $O_{b}$ called the \textbf{blue orthocenter}. The
three red altitudes $a_{1}^{r},a_{2}^{r}$ and $a_{3}^{r}$ meet in a point
$O_{r}$ called the \textbf{red orthocenter}. The three green altitudes
$a_{1}^{g},a_{2}^{g}$ and $a_{3}^{g}$ meet in a point $O_{g}$ called the
\textbf{green orthocenter}. For $A_{1}\equiv\left[  x_{1},y_{1}\right]  $,
$A_{2}\equiv\left[  x_{2},y_{2}\right]  $ and $A_{3}\equiv\left[  x_{3}%
,y_{3}\right]  $ these points are given by
\begin{align*}
O_{b}  & =\left[  \frac{\left[  x_{1}x_{2}y_{2}\right]  ^{-}+\left[
y_{1}y_{2}^{2}\right]  ^{-}}{\left[  x_{1}y_{2}\right]  ^{-}},\frac{\left[
x_{1}y_{1}y_{2}\right]  ^{-}+\left[  x_{1}^{2}x_{2}\right]  ^{-}}{\left[
x_{1}y_{2}\right]  ^{-}}\right] \\
O_{r}  & =\left[  \frac{\left[  x_{1}x_{2}y_{2}\right]  ^{-}-\left[
y_{1}y_{2}^{2}\right]  ^{-}}{\left[  x_{1}y_{2}\right]  ^{-}},\frac{\left[
x_{1}y_{1}y_{2}\right]  ^{-}-\left[  x_{1}^{2}x_{2}\right]  ^{-}}{\left[
x_{1}y_{2}\right]  ^{-}}\right] \\
O_{g}  & =\left[  \frac{\left[  x_{1}^{2}y_{2}\right]  ^{-}+\left[  x_{1}%
x_{2}y_{1}\right]  ^{-}}{\left[  x_{1}y_{2}\right]  ^{-}},\frac{\left[
x_{1}y_{2}^{2}\right]  ^{-}-\left[  x_{1}y_{1}y_{2}\right]  ^{-}}{\left[
x_{1}y_{2}\right]  ^{-}}\right]  .%
%TCIMACRO{\TeXButton{Proof box}{{\hspace{.1in} \rule{0.5em}{0.5em}}}}%
%BeginExpansion
{\hspace{.1in} \rule{0.5em}{0.5em}}%
%EndExpansion
\end{align*}

\end{theorem}

\section*{Circumcenters}

When $A_{1}$ and $A_{2}$ are distinct points with $l=A_{1}A_{2},$ and $A$ is
the midpoint of $A_{1}$ and $A_{2}$, then the blue, red and green altitudes
from $A$ to $l$ are respectively called the \textbf{blue},\textbf{\ red
}and\textbf{\ green perpendicular bisectors }of the side $\overline{A_{1}%
A_{2}}$.

\begin{theorem}
[Perpendicular bisectors]If $A_{1}\equiv\left[  x_{1},y_{1}\right]  $ and
$A_{2}\equiv\left[  x_{2},y_{2}\right]  $ are distinct points then the blue,
red and green perpendicular bisectors of $\overline{A_{1}A_{2}}$ have
respective equations%
\begin{align*}
\left(  x_{1}-x_{2}\right)  x+\left(  y_{1}-y_{2}\right)  y  & =\frac
{x_{1}^{2}-x_{2}^{2}+y_{1}^{2}-y_{2}^{2}}{2}\\
\left(  x_{1}-x_{2}\right)  x-\left(  y_{1}-y_{2}\right)  y  & =\frac
{x_{1}^{2}-x_{2}^{2}-y_{1}^{2}+y_{2}^{2}}{2}\\
\left(  y_{2}-y_{1}\right)  x+\left(  x_{2}-x_{1}\right)  y  & =y_{2}%
x_{2}-x_{1}y_{1}.%
%TCIMACRO{\TeXButton{Proof box}{{\hspace{.1in} \rule{0.5em}{0.5em}}}}%
%BeginExpansion
{\hspace{.1in} \rule{0.5em}{0.5em}}%
%EndExpansion
\end{align*}

\end{theorem}

Given a triangle $\overline{A_{1}A_{2}A_{3}}$ we may construct the blue, red
and green perpendicular bisectors of the three sides, denoted by $b_{m}%
^{b},b_{m}^{r}$ and $b_{m}^{g}$ respectively for $m=1,2$ and $3,$ where
$b_{1}^{b}$ for example is the blue perpendicular bisector of the side
$\overline{A_{2}A_{3}}$ and so on.$\allowbreak$

\begin{theorem}
[Circumcenter formulas]The three blue perpendicular bisectors $b_{1}^{b}%
,b_{2}^{b}$ and $b_{3}^{b}$ meet in a point $C_{b}$ called the \textbf{blue
circumcenter}. The three red perpendicular bisectors $b_{1}^{r},b_{2}^{r}$ and
$b_{3}^{r}$ meet in a point $C_{r}$ called the \textbf{red circumcenter}. The
three green perpendicular bisectors $b_{1}^{g},b_{2}^{g}$ and $b_{3}^{g}$ meet
in a point $C_{g}$ called the \textbf{green circumcenter}. For $A_{1}%
\equiv\left[  x_{1},y_{1}\right]  $, $A_{2}\equiv\left[  x_{2},y_{2}\right]  $
and $A_{3}\equiv\left[  x_{3},y_{3}\right]  $, these points are given by%
\begin{align*}
C_{b}  & =\left[  \frac{\left[  x_{1}^{2}y_{2}\right]  ^{-}+\left[  y_{1}%
^{2}y_{2}\right]  ^{-}}{2\left[  x_{1}y_{2}\right]  ^{-}},\frac{\left[
x_{1}y_{2}^{2}\right]  ^{-}+\left[  x_{1}x_{2}^{2}\right]  ^{-}}{2\left[
x_{1}y_{2}\right]  ^{-}}\right] \\
C_{r}  & =\left[  \frac{\left[  x_{1}^{2}y_{2}\right]  ^{-}-\left[  y_{1}%
^{2}y_{2}\right]  ^{-}}{2\left[  x_{1}y_{2}\right]  ^{-}},\frac{\left[
x_{1}y_{2}^{2}\right]  ^{-}-\left[  x_{1}x_{2}^{2}\right]  ^{-}}{2\left[
x_{1}y_{2}\right]  ^{-}}\right] \\
C_{g}  & =\left[  \frac{\left[  x_{1}x_{2}y_{2}\right]  ^{-}}{\left[
x_{1}y_{2}\right]  ^{-}},\frac{\left[  x_{1}y_{1}y_{2}\right]  ^{-}}{\left[
x_{1}y_{2}\right]  ^{-}}\right]  .%
%TCIMACRO{\TeXButton{Proof box}{{\hspace{.1in} \rule{0.5em}{0.5em}}}}%
%BeginExpansion
{\hspace{.1in} \rule{0.5em}{0.5em}}%
%EndExpansion
\end{align*}

\end{theorem}

\begin{theorem}
[Circumcenters as midpoints]For any triangle a coloured circumcenter is the
midpoint of the two orthocenters of the other two colours.
\end{theorem}

\begin{proof}
This follows from the Orthocenter formulas and Circumcenter formulas.
\end{proof}

\section*{Nine-point centres}

Suppose that the respective midpoints of a triangle $\overline{A_{1}A_{2}%
A_{3}}$ are $M_{m}$ for $m=1,2$ and $3$, where $M_{1}$ is the midpoint of the
side $\overline{A_{2}A_{3}}$ and so on. We let $N_{b},N_{r}$ and $N_{g}$ be
the blue, red and green circumcenters respectively of the triangle
$\overline{M_{1}M_{2}M_{3}}$, and call these the \textbf{blue},\textbf{\ red
}and\textbf{\ green nine-point centers} of the original triangle
$\overline{A_{1}A_{2}A_{3}}$.

\begin{theorem}
[Nine-point center formulas]For $A_{1}\equiv\left[  x_{1},y_{1}\right]  $,
$A_{2}\equiv\left[  x_{2},y_{2}\right]  $ and $A_{3}\equiv\left[  x_{3}%
,y_{3}\right]  $, the blue, red and green nine-point centers of $\overline
{A_{1}A_{2}A_{3}}$ are%
\begin{align*}
N_{b}  & =\left[  \frac{\left[  x_{1}^{2}y_{2}\right]  ^{-}-\left[  y_{1}%
^{2}y_{2}\right]  ^{-}+2\left[  x_{1}x_{2}y_{2}\right]  ^{-}}{4\left[
x_{1}y_{2}\right]  ^{-}},\frac{\left[  x_{1}y_{2}^{2}\right]  ^{-}-\left[
x_{1}x_{2}^{2}\right]  ^{-}+2\left[  x_{1}y_{1}y_{2}\right]  ^{-}}{4\left[
x_{1}y_{2}\right]  ^{-}}\right] \\
N_{r}  & =\left[  \frac{\left[  x_{1}^{2}y_{2}\right]  ^{-}+\left[  y_{1}%
^{2}y_{2}\right]  ^{-}+2\left[  x_{1}x_{2}y_{2}\right]  ^{-}}{4\left[
x_{1}y_{2}\right]  ^{-}},\frac{\left[  x_{1}y_{2}^{2}\right]  ^{-}+\left[
x_{1}x_{2}^{2}\right]  ^{-}+2\left[  x_{1}y_{1}y_{2}\right]  ^{-}}{4\left[
x_{1}y_{2}\right]  ^{-}}\right] \\
N_{g}  & =\left[  \frac{\left[  x_{1}^{2}y_{2}\right]  ^{-}}{2\left[
x_{1}y_{2}\right]  ^{-}},\frac{\left[  x_{1}y_{2}^{2}\right]  ^{-}}{2\left[
x_{1}y_{2}\right]  ^{-}}\right]  .%
%TCIMACRO{\TeXButton{Proof box}{{\hspace{.1in} \rule{0.5em}{0.5em}}}}%
%BeginExpansion
{\hspace{.1in} \rule{0.5em}{0.5em}}%
%EndExpansion
\end{align*}

\end{theorem}

\begin{theorem}
[Nine-point centers as midpoints]In any triangle a coloured nine-point center
is the midpoint of the two circumcenters of the other two colours.
\end{theorem}

\begin{proof}
This follows from the Circumcenter formulas and Nine-point center formulas.
\end{proof}

\section*{The $\Omega$-triangle and the Euler lines}

The $\Omega$\textbf{-triangle }of\textbf{\ }a triangle $\overline{A_{1}%
A_{2}A_{3}}$ is the triangle $\Omega\equiv\Omega\left(  \overline{A_{1}%
A_{2}A_{3}}\right)  \equiv\overline{O_{b}O_{r}O_{g}}$ of orthocenters of
$\overline{A_{1}A_{2}A_{3}}$. From the theorems of the last two sections, the
corresponding midpoints of the sides of $\Omega$ are $C_{b},C_{r}$ and
$C_{g},$ with $C_{b}$ the midpoint of $O_{r}$ and $O_{g}$ etc., and the
midpoints of the triangle $\overline{C_{b}C_{r}C_{g}}$ are $N_{b},N_{r}$ and
$N_{g},$ with $N_{b}$ the midpoint of $C_{r}$ and $C_{g}$ etc. We also know
that the centroid of $\Omega$ is the same as the centroid $G$ of the original
triangle $\overline{A_{1}A_{2}A_{3}}$.

\begin{theorem}
[Blue Euler line]The points $O_{b},N_{b},G$ and $C_{b}$ lie on a line called
the \textbf{blue Euler line}. Furthermore $N_{b}$ is the midpoint of $O_{b}$
and $C_{b},$ and we have the affine relations%
\[
G=\frac{1}{3}O_{b}+\frac{2}{3}C_{b}=\frac{1}{3}C_{b}+\frac{2}{3}N_{b}.
\]

\end{theorem}

\begin{proof}
This follows from the Orthocenter, Circumcenter and Nine-point center formulas.
\end{proof}

\begin{theorem}
[Red Euler line]The points $O_{r},N_{r},G$ and $C_{r}$ lie on a line called
the \textbf{red Euler line}. Furthermore $N_{r}$ is the midpoint of $O_{r}$
and $C_{r},$ and%
\[
G=\frac{1}{3}O_{r}+\frac{2}{3}C_{r}=\frac{1}{3}C_{r}+\frac{2}{3}N_{r}.
\]

\end{theorem}

\begin{proof}
Likewise.
\end{proof}

\begin{theorem}
[Green Euler line]The points $O_{g},N_{g},G$ and $C_{g}$ lie on a line called
the \textbf{green Euler line}. Furthermore $N_{g}$ is the midpoint of $O_{g}$
and $C_{g}, $ and%
\[
G=\frac{1}{3}O_{g}+\frac{2}{3}C_{g}=\frac{1}{3}C_{g}+\frac{2}{3}N_{g}.
\]

\end{theorem}

\begin{proof}
Likewise.
\end{proof}

The geometry of the $\Omega$-triangle clarifies the various ratios occurring
along points on the Euler lines, since these are just the medians of $\Omega$.
The lines joining the circumcenters are the lines of the medial triangle of
$\Omega,$ and so are parallel to the lines of $\Omega$.

\section*{Circles}

A \textbf{blue, red or green circle} is an equation $c$ in $x$ and $y$ of the
form
\begin{align*}
\left(  x-x_{0}\right)  ^{2}+\left(  y-y_{0}\right)  ^{2}  & =K\\
\left(  x-x_{0}\right)  ^{2}-\left(  y-y_{0}\right)  ^{2}  & =K\\
2\left(  x-x_{0}\right)  \left(  y-y_{0}\right)   & =K
\end{align*}
respectively, where the point $\left[  x_{0},y_{0}\right]  $ is then unique
and called the \textbf{centre} of $c$, and $K$ is the \textbf{quadrance} of
$c$. A blue circle is an ordinary Euclidean circle. Red and green circles are
more usually described as rectangular hyperbolas. A red circle has asymptotes
parallel to the lines with equations $y=\pm x,$ and a green circle has
asymptotes parallel to the coordinate axes.

\begin{theorem}
[Circumcircles]If $A_{1},A_{2}$ and $A_{3}$ are three distinct non-collinear
points, then there are unique blue, red and green circles passing through
$A_{1},A_{2}$ and $A_{3}$. $%
%TCIMACRO{\TeXButton{Proof box}{{\hspace{.1in} \rule{0.5em}{0.5em}}}}%
%BeginExpansion
{\hspace{.1in} \rule{0.5em}{0.5em}}%
%EndExpansion
$
\end{theorem}

The circles above will be called respectively the \textbf{blue},\textbf{\ red
}and\textbf{\ green circumcircles} of the triangle $\overline{A_{1}A_{2}A_{3}%
},$ while the circumcircles of the triangle of midpoints $M_{1},M_{2},M_{3}$
of the triangle $\overline{A_{1}A_{2}A_{3}}$ will be called respectively the
\textbf{blue},\textbf{\ red }and \textbf{green nine-point circles }of the
triangle $\overline{A_{1}A_{2}A_{3}}$.

\begin{theorem}
[Orthocenters on circumcircles]Any coloured orthocenter of a triangle
$\overline{A_{1}A_{2}A_{3}}$ lies on the circumcircles of the other two
colours. $%
%TCIMACRO{\TeXButton{Proof box}{{\hspace{.1in} \rule{0.5em}{0.5em}}}}%
%BeginExpansion
{\hspace{.1in} \rule{0.5em}{0.5em}}%
%EndExpansion
$
\end{theorem}

\begin{theorem}
[Nine-point circles]Any coloured nine-point circle of a triangle
$\overline{A_{1}A_{2}A_{3}}$ passes through the feet of the altitudes of that
colour, as well as the midpoints of the segments from the same coloured
orthocenter to the points $A_{1},A_{2}$ and $A_{3}.$ In addition it passes
through the circumcenters of $\overline{A_{1}A_{2}A_{3}}$ of the other two
colours. $%
%TCIMACRO{\TeXButton{Proof box}{{\hspace{.1in} \rule{0.5em}{0.5em}}}}%
%BeginExpansion
{\hspace{.1in} \rule{0.5em}{0.5em}}%
%EndExpansion
$
\end{theorem}

The following figure shows some of the other points on the nine-point circles
of different colours. Others are off the page.%

%TCIMACRO{\FRAME{dhF}{5.0959in}{3.716in}{0pt}{}{}{chromoninepoint.eps}%
%width 5.0959in;  height 3.716in;  depth 0pt;  original-width 7.0673in;
%original-height 5.144in;  cropleft "0";  croptop "1";  cropright "1";
%cropbottom "0";  filename 'ChromoNinePoint.eps';file-properties "XNPEU";}}}%
%BeginExpansion
\begin{center}
\includegraphics[
height=3.716in,
width=5.0959in
]%
{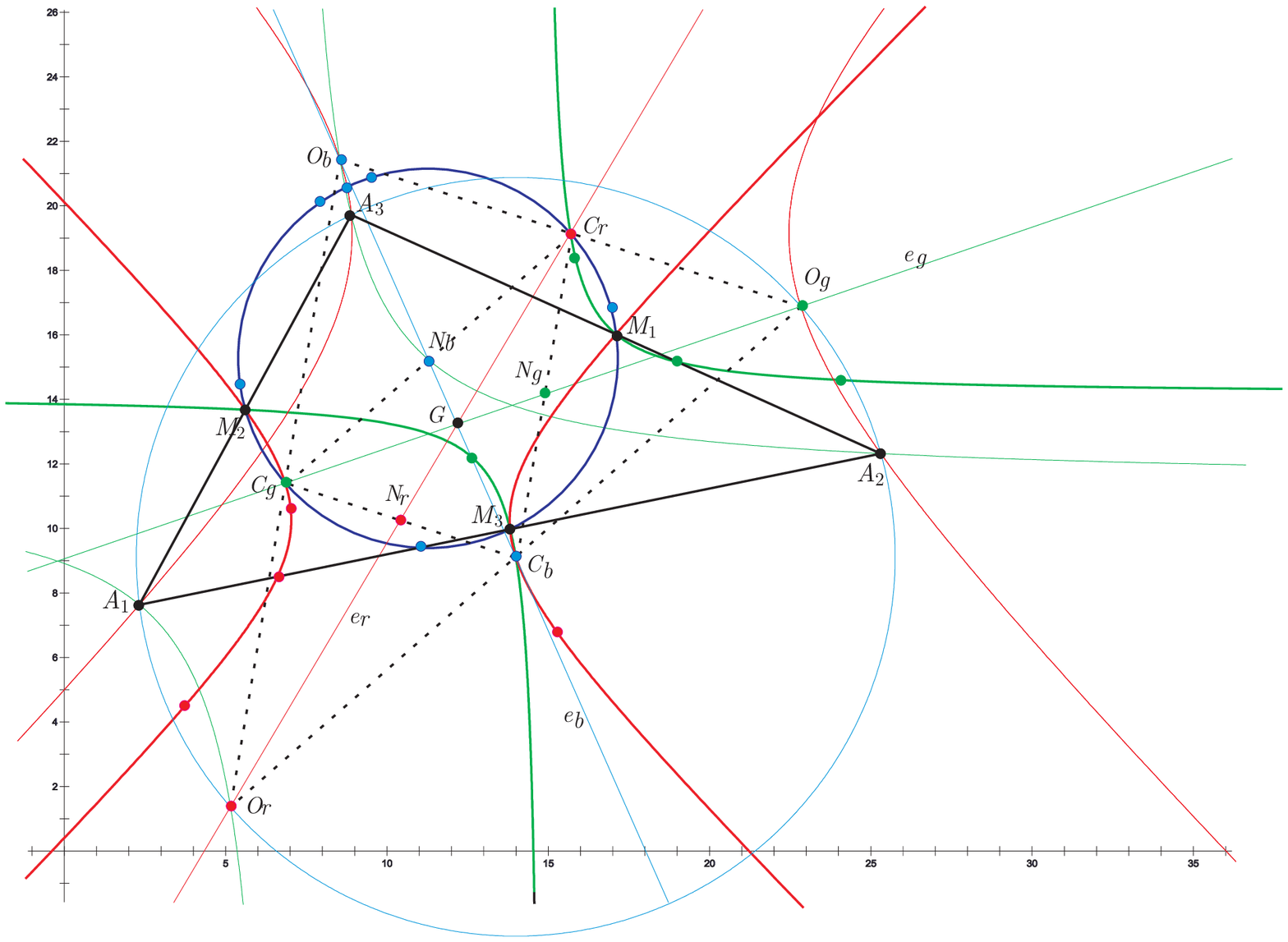}%
\end{center}
%EndExpansion

Hopefully this taste of chromogeometry will encourage others to explore this
rich new realm. See \cite{Wildberger3} and \cite{Wildberger 4} for more in
this direction.

\end{document}